\documentclass[11pt, english]{amsart}
\usepackage{babel,amssymb,amsmath,latexsym,amsthm,amscd,eucal,enumerate,verbatim}  
\usepackage[dvips]{graphicx}
\usepackage[all]{xy}
\usepackage[left=2.5cm,top=2.5cm, bottom=2.5cm,right=2.5cm]{geometry}
\usepackage{setspace}

\newtheorem{thm}{Theorem}

\newtheorem{lem}[thm]{Lemma}
\newtheorem{cor}[thm]{Corollary}

\newcommand\bp{\noindent{\it Proof.}\ }

\title{Extension of positive maps}

\author{Erling St{\o}rmer}

\date{14-3-2018}
\begin{document}
\maketitle

\begin{abstract}

{We prove two extension theorems for positive maps from operator systems into matrix algebras}

\end{abstract}

In the theory of positive maps of operator algebras Arveson`s Extension Theorem for completely positive maps plays a major role.  In the paper \cite {S} the author extended his result to maps with more general positivity properties, the main theorem being included in Theorem 5.2.3 in the book \cite{St}. However, it was pointed out by D. Chruscinski to the author that V. Paulsen had a counter-example for general positive maps to this theorem in his book \cite{P}, Example 2.2, see also \cite{B}.

In the present note we show two corrected versions of the above Theorem 5.2.3.  The first theorem is very close to the original except that we restrict the operator system to consist only of self-adjoint matrices.  This is due to the fact that Krein`s Extension Theorem \cite{St}, Theorem A.3.1, is formulated for real spaces, while in Theorem 5.2.3 we wrongly applied it to complex spaces. 

The other extension theorem is an extension theorem for general positive maps and is closely related to the fact that a unital linear map $\phi$ of a C*-algebra into $B(H)$ is positive if and only if $\parallel \phi\parallel = \parallel \phi(1)\parallel = 1.$

Our basic reference for these notes is the book \cite{St}.  We recall some concepts which we shall use.  An operator system is a complex self-adjoint linear subspace  A of operators in $B(H)$ such that $1\in A$.  A mapping cone  $\mathcal{C}$ on $H$ is a closed convex subcone of the cone of positive maps of $B(H)$ into itself such that $\phi \in \mathcal{C}$ implies $ \alpha \circ \phi \circ \beta \in \mathcal{C}$ for all completely positive maps $\alpha, \beta$ of $B(H)$ into itself.  For simplicity we assume $H$ and $K$ are finite dimensional Hilbert spaces, and we let $Tr$ denote usual trace on $B(K)$ or on $B(H) \otimes B(K)$ if there is no confusion. Let $ \phi $ be a linear map of $A$ into $B(K).$ Then the dual functional $\tilde \phi$  on $B(H)\otimes B(K)$ is defined by the formula
$$
\tilde \phi(a \otimes b) = Tr(\phi(a) b^t) = Tr(C_{\phi} a\otimes b),
$$
where $b^t$ is the transpose of $b$, and $ C_{\phi}$ is the Choi matrix for $\phi$, see \cite{St}, Definition 4.1.1. If $\mathcal{C}$ is a mapping cone as above, we denote by $P(A,\mathcal{C})$.  the cone
$$
P(A,\mathcal{C}) = \{ x \in (A \otimes B(K))_{sa} : \iota \otimes \alpha(x) \geq 0, \forall \alpha \in \mathcal{C} \},
$$
where $\iota$ is the identity map on $B(H)$. We say $\phi$ is $\mathcal{C}$- positive if $\tilde \phi$ is positive on $P(A,\mathcal{C}).$

\section{Main results}

As mentioned above our first extension theorem is for the real version of operator systems .  If $A$ is a real linear subspace of the self-adjoint  operators in $B(H)$ containing the identity $1$, we say $A$ is a real operator system. The definition of $\mathcal{C}$-positive maps still make sense for real operator systems. $H$ and $K$ are still finite dimensional Hilbert spaces.

 \begin{thm}\label {thm}
 Let $A$ be a real operator system contained in $B(H)$.  Let $\phi$ be a $\mathcal{C}$ positive map of $A$ into $B(K)$ for a mapping cone $\mathcal{C}$.  Then there exists
 a $\mathcal{C}$ positive map $\psi$ of $B(H)$ into $B(K)$ such that $\psi(a) = \phi(a)$ for $a \in A$.
 \end{thm}
 \bp 
 Let 
 $P = P(B(H), \mathcal{C}) $ be defined as above.  Then $P(A,\mathcal{C}) = P \cap (A \otimes B(K))_{sa}$.  Since $\phi$ is $\mathcal{C}$-positive  its dual functional $\tilde \phi$ is positive on $P\cap (A \otimes B(K))_{sa}$.  By \cite{St}, Lemma 5.2.1, $1\otimes 1$ is an interior point of $P$. Thus by Kreins Extension Theorem  \cite{St}, Theorem A.3.1, $\tilde \phi$ has an extension to a real linear functional $\tilde \psi _o$ on $(B(H) \otimes B(K))_{sa},$ which is positive on $P$.  Define a complex linear functional $\tilde \psi$ on $B(H) \otimes B(K)$ by 
  $$
 \tilde\psi (a + ib) = \tilde\psi_{o} a) + i \tilde\psi_{o} (b),\ \forall \ a,b \in (B(H)\otimes. B(K))_{sa}.
  $$
 A straightforward computation shows that 
 $$
 \tilde\psi(\lambda x) = \lambda \tilde\psi(x),\  \forall \ x\in B(H)\otimes B(K).
 $$
 Thus $\tilde \psi$ is a complex linear functional on $B(H)\otimes B(K)$ which is positive on $P$ .  But then there exists an operator $C\in B(H) \otimes B(K)$ such that 
 $$
 \tilde\psi(X) = Tr(Cx) \ \forall \ x\in B(H)\otimes B(K).
 $$
 By \cite{St},Lemmas 4.22 and 4.23 there exists a linear map $\psi$ of $B(H)$ into $B(K) $ such that $C = C_{\psi}$ is the Choi matrix for $\psi$. Then 
 $$
 Tr(C_{\psi}x) \geq 0 \forall x\in P.
 $$
 Thus $\psi$ is $\mathcal{C}$-positive, and if $a\in A, b\in B(K)$ then
 \begin{eqnarray*}
 Tr(\psi(a)b^t) &=& Tr(C_{\psi} a\otimes b)\\
 &=&\tilde\psi(a\otimes b)\\
 &=& \tilde\phi(a\otimes b)\\
 &=& Tr(\phi(a)b^t).
  \end{eqnarray*}
 Since this holds for all $b\in B(K)$, $\psi(a) = \phi(a)$ completing the proof.
\medskip
 
 We need $H$ and $K$ finite dimensional in order to choose the operator $C$ in the proof.  The theorem can be extended to the case when $K$ infinite dimensional by the same proof as that of \cite{St}, Theorem 5.2.3. 
 
 \begin{cor}\label{cor}
 Let $B$ be a C*-subalgebra of $B(H)$ and $\phi$ a $\mathcal{C}$-positive map of $B$ into $B(K)$.  Then $\phi$ has a $\mathcal{C}$-positive extension $\psi$ of $B(H)$ into $B(K)$.
 \end{cor}
 \bp  Let $A = B_{sa}$.  By the theorem applied to $A$ the corollary follows.

 In our next theorem we use the Hahn-Banach Theorem for extending $\tilde \phi $ to $\tilde \psi$ as above.  But we have to restrict attention to the mapping cone consisting of all positive maps.  The proof is simplest when we assume the map is unital.  The next lemma shows that we can reduce to this case.  Recall that $Ad V$ is the map $ a \longrightarrow V^*aV$.
 
 \begin {lem} \label {lem }
 Let $A$ be an operator system $A \subset B(H)$ and $\phi$ a positive map of $A$ into $B(K)$, where $H$ and $K$ are finite dimensional Hilbert spaces.  Let $V$ be an invertible operator in $B(K)$.  Then $\phi$ has an extension to a positive map $\psi$  of $B(H)$ into $B(K)$ if and only if  $Ad V \circ \phi$ has a positive extension.
 \end {lem}
 
 \bp. If $\psi$ is an extension of $\phi$ then $Ad V\circ \psi$ is an extension of $AdV\circ \phi$.  Conversely, if $\eta$ is an extension of $AdV \circ \phi$, then $\psi = Ad V^{-1} \circ \eta$ is an extension of $\phi$.
 \medskip
 
 If $\phi$ is a positive map of $A$ into $B(K)$ let $p$ be the range projection of $\phi(1)$.  Considering $Ad p \circ \phi$ mapping $A$ into $pB(K)p$ instead of $\phi$, we may assume $\phi(1)$ is strictly positive, so invertible.
 
 \begin {thm}\label {thm}
 Let $H$ and $K$ be finite dimensional Hilbert spaces. Let $A\subset B(H)$ be an operator system, and let $\phi$ be a positive map of $A$ into $B(K)$ with $\phi(1)$ invertible.  Then $\phi$ has a positive extension $\psi$ of $B(H)$ into $B(K)$ if and only if $\parallel Ad \phi(1)^{-1/2} \circ \phi \parallel = 1.$
 \end{thm}
 
 \bp  The map $Ad \phi(1)^{-1/2} \circ \phi$ is unital, so by Lemma 3 we may replace $\phi$ by $Ad \phi(1)^{-1/2} \circ \phi$, and thus assume $\phi(1) = 1$.  
 
 Suppose $\psi$ is a positive extension of $\phi$ to a positive map of $B(H)$ into $B(K)$.  Since $B(H)$ is a C*-algebra, and $\psi$ is positive, by \cite{St} Theorem 1.3.3
 
$$ \parallel \psi\parallel = \parallel \psi(1)\parallel =\parallel \phi(1) \parallel \leq \parallel \phi\parallel \leq \parallel \psi \parallel.$$
 
 Thus $\parallel \phi\parallel = \parallel \phi(1) \parallel = 1.$
 
 Conversely, assume $\parallel \phi\parallel = \parallel \phi(1) \parallel =1$.  By \cite{St}, Lemma 4.2.2 the  map $\phi \longrightarrow \tilde \phi$ is an isometry where $B(H) \otimes B(K)$ has the projective norm , $\parallel a \otimes b\parallel = \parallel a\parallel \parallel b \parallel_1$, where $\parallel b\parallel_1$ is the trace norm of $b$.  By the norm version of the Hahn-Banach Theorem \cite{KR} Theorem 1.6.1, $\tilde \phi$ has an extension from $A\otimes B(K)$ to a linear functional $\tilde \psi$ on $B(H) \otimes B(K)$ such that $\parallel \tilde \psi\parallel = \parallel \tilde \phi\parallel$, which is equal to $\parallel \phi\parallel = 1$ by \cite{St} Lemma 4.2.2.  By the same lemma $\tilde \phi$ is positive on $A^+ \otimes B(K)^+$ hence so is the extension $\tilde\psi$.  Thus for $b\in B(K)^+$ the linear functional on $B(H)$,
 $$
 f_{b}(a)  = \tilde \psi(a \otimes b)
 $$ 
 is positive on $A$ and is an extension of the functional $a\longrightarrow \tilde \phi(a\otimes b)$.  Furthermore, for $a \in B(H)$
 \begin{eqnarray*}
 \mid f_b(a) \mid &=& \mid \tilde\psi(a \otimes b)\mid \\
 &\leq& \parallel \tilde\psi\parallel  \parallel a\otimes b\parallel \\
 &=& \parallel \tilde \phi\parallel \parallel a \parallel \parallel b\parallel_1\\
 &=& \parallel a\parallel \parallel b \parallel_1
  \end{eqnarray*}
 Thus $f_b$ satisfies, using that $\phi(1) = 1$,
 \begin{eqnarray*}
 \parallel b \parallel _1&\geq& \parallel f_b \parallel\\
 &=& sup \{ \mid \tilde\psi(a\otimes b)\mid : \parallel a \parallel \leq 1, a \in B(H) \}\\
 &\geq& \mid \tilde\psi(1\otimes b)\mid\\
 &=& \tilde\phi(1\otimes b)\\
 &=& Tr(\phi(1)b^t)\\
 &=& Tr(b^t)\\
 &=& \parallel b\parallel_1.
 \end{eqnarray*}
 It follows that 
 $$
 \parallel f_b\parallel = \parallel b\parallel_1 = f_b(1).
 $$
 Thus by \cite{KR}, Theorem 4.3.2, $f_b \geq 0$ 0n $B(H),$ and so for $a\geq 0, b\geq 0$
 $$
 0 \leq f_b(a) = \tilde \psi(a \otimes b).
 $$
 Thus by \cite{St}, Lemma 4.2.2, $\tilde\psi$ is positive on $B(H)^+ \otimes B(K)^+$.  As in the proof of Theorem 1 there is a positive map $\psi$ of $B(H)$ into$B(K)$ such that 
 $$
 \tilde\psi(a\otimes b) = Tr(C_{\psi} a\otimes b) = Tr(\psi(a) b^t).
 $$
 If $a \in A$ and $b \in B(K)$ it follows that
 $$
 Tr(\psi(a) b^t) = \tilde\psi(a\otimes b) = \tilde\phi(a\otimes b) = Tr(\phi(a)b^t),
 $$
 Hence $\psi(a) =\phi(a)$, and $\psi$ is the desired extension of $\phi$.  The proof is complete.
 
 The above proof yields a proof of the Arvesons Extension Theorem in the finite dimensional ase, because if $\phi$ is completely positive then by \cite{St}, Theorem 4.2.7, $\tilde \phi \geq 0,$ hence by \cite{KR},Theorem  4.3.13 we can choose $\tilde \psi $ to be positive on $B(H)\otimes B(K)$.  Thus $C_{\psi} \geq 0$, and so $\psi$ is completely positive by \cite{St}, Theorem 4.1.8.
 
 By the same techniques as in the proof of \cite{St}, Theorem 5.3.2, we can extend the last theorem to the case when K is infinite dimensional.  The case when H is infinite dimensional is more complex because we need $\tilde \psi$ to be ultraweakly continuous in order to choose $C_{\psi} $as a Radon-Nikodym derivative for $\tilde\psi.$
 \medskip
 
 Acknowledgements
 The author wants to thank D. Chrusciniski for useful correspondence on the subject and E. Alfsen for very useful discussions on the proofs of the theorems in the paper.

Department of Mathematics,  University of Oslo, 0316 Oslo, Norway.

e-mail  erlings@math.uio.no

\end{document}